\documentstyle[twoside,11pt]{article}
\begin{document}
\pagestyle{myheadings}
\markboth{Gantz and  Steer}{Stable
 Parabolic Bundles over Elliptic Surfaces}

\title{Stable Parabolic Bundles over Elliptic Surfaces
      and over Orbifold Riemann Surfaces}
\author{Christian Gantz$^{1}$, Brian Steer \\
         Mathematical Institute, Oxford OX1 3LB, UK \\
       (gantz@maths.ox.ac.uk) }

\maketitle

\newcommand{\exactss}[5]
{\mbox{$0 \rightarrow #1 \stackrel{#4}{\rightarrow} #2
\stackrel{#5}{\rightarrow} #3 \rightarrow 0$}}

\newcommand{\fourhorbbb}[4]
{\mbox{$#1 \rightarrow #2 \rightarrow #3 \rightarrow #4$}}

\newcommand{\fourhorbll}[6]
{\mbox{$#1 \rightarrow #2 \stackrel{#5}{\longrightarrow} #3
\stackrel{#6}{\longrightarrow}
 #4$}}

\newcommand{\fourhorsss}[7]
{\mbox{$#1 \stackrel{#5}{\rightarrow} #2
\stackrel{#6}{\rightarrow} #3 \stackrel{#7}{\rightarrow}
 #4$}}

\newcommand{\lifta}[6]
{\setlength{\unitlength}{1.0mm}
\begin{center}
\begin{picture}(14,56)(-7,-28)

\putcc{13}{13}{#1}
\putcc{-13}{-13}{#2}
\putcc{13}{-13}{#3}

\multiput(-8,-8)(5,5){3}{\line(1,1){4}}
\put(8,8){\vector(1,1){0}}
\putbr{-1}{-1}{#4}
\put(-6,-13){\vector(1,0){12}}
\puttc{0}{-9}{#6}
\put(13,8){\vector(0,-1){16}}
\putcl{15}{0}{#5}
\end{picture}
\end{center}}

\newcommand{\quadrata}[8]
{\setlength{\unitlength}{1.0mm}
\begin{center}
\begin{picture}(14,46)(-7,-23)

\putcc{-13}{13}{#1}
\putcc{13}{13}{#2}
\putcc{-13}{-13}{#3}
\putcc{13}{-13}{#4}

\put(-6,13){\vector(1,0){12}}
\putbc{0}{15}{#5}
\put(-6,-13){\vector(1,0){12}}
\putbc{0}{-11}{#8}
\put(-13,8){\vector(0,-1){16}}
\putcr{-15}{0}{#6}
\put(13,8){\vector(0,-1){16}}
\putcl{15}{0}{#7}
\end{picture}
\end{center}}

\newcommand{\quadratb}[8]
{\setlength{\unitlength}{1.0mm}
\begin{center}
\begin{picture}(14,56)(-7,-28)

\putcc{-18}{13}{#1}
\putcc{18}{13}{#2}
\putcc{-18}{-13}{#3}
\putcc{18}{-13}{#4}

\put(-9,13){\vector(1,0){18}}
\putbc{0}{15}{#5}
\put(-9,-13){\vector(1,0){18}}
\putbc{0}{-11}{#8}
\put(-18,8){\vector(0,-1){16}}
\putcr{-20}{0}{#6}
\put(18,8){\vector(0,-1){16}}
\putcl{20}{0}{#7}
\end{picture}
\end{center}}

\newcommand{\quadratc}[8]
{\setlength{\unitlength}{1.0mm}
\begin{center}
\begin{picture}(14,56)(-7,-28)

\putcc{-23}{13}{#1}
\putcc{23}{13}{#2}
\putcc{-23}{-13}{#3}
\putcc{23}{-13}{#4}

\put(-11,13){\vector(1,0){22}}
\putbc{0}{15}{#5}
\put(-11,-13){\vector(1,0){22}}
\putbc{0}{-11}{#8}
\put(-23,8){\vector(0,-1){16}}
\putcr{-25}{0}{#6}
\put(23,8){\vector(0,-1){16}}
\putcl{25}{0}{#7}
\end{picture}
\end{center}}

\newcommand{\quadratd}[8]
{\setlength{\unitlength}{1.0mm}
\begin{center}
\begin{picture}(14,56)(-7,-28)

\putcc{-29}{13}{#1}
\putcc{29}{13}{#2}
\putcc{-23}{-13}{#3}
\putcc{23}{-13}{#4}

\put(-11,13){\vector(1,0){22}}
\putbc{0}{15}{#5}
\put(-11,-13){\vector(1,0){22}}
\putbc{0}{-11}{#8}
\put(-23,8){\vector(0,-1){16}}
\putcr{-25}{0}{#6}
\put(23,8){\vector(0,-1){16}}
\putcl{25}{0}{#7}
\end{picture}
\end{center}}

\newcommand{\quadrate}[8]
{\setlength{\unitlength}{1.0mm}
\begin{center}
\begin{picture}(14,56)(-7,-28)

\putcc{-25}{13}{#1}
\putcc{29}{13}{#2}
\putcc{-25}{-13}{#3}
\putcc{29}{-13}{#4}

\put(-11,13){\vector(1,0){22}}
\putbc{0}{15}{#5}
\put(-11,-13){\vector(1,0){22}}
\putbc{0}{-11}{#8}
\put(-23,8){\vector(0,-1){16}}
\putcr{-25}{0}{#6}
\put(23,8){\vector(0,-1){16}}
\putcl{25}{0}{#7}
\end{picture}
\end{center}}

\newcommand{\quadratf}[8]
{\setlength{\unitlength}{1.0mm}
\begin{center}
\begin{picture}(14,56)(-7,-28)

\putcc{-35}{13}{#1}
\putcc{34}{13}{#2}
\putcc{-35}{-13}{#3}
\putcc{34}{-13}{#4}

\put(-11,13){\vector(1,0){22}}
\putbc{0}{15}{#5}
\put(-11,-13){\vector(1,0){22}}
\putbc{0}{-11}{#8}
\put(-23,8){\vector(0,-1){16}}
\putcr{-25}{0}{#6}
\put(23,8){\vector(0,-1){16}}
\putcl{25}{0}{#7}
\end{picture}
\end{center}}

\newcommand{\threeverta}[5]
{\begin{picture}(140,140)(-175,-70)
\setlength{\unitlength}{0.7pt}
\put(-50,50){$#1$}
\put(-50,-50){$#2$}
\put(-50,-150){$#3$}
\put(-42,37){\vector(0,-1){60}}
\put(-73,0){$#4$}
\put(-42,-63){\vector(0,-1){60}}
\put(65,0){$#5$}
\end{picture}}

\newcommand{\threehorbl}[4]
{\mbox{$#1 \rightarrow #2 \stackrel{#4}{\longrightarrow} #3$}}

\newcommand{\threehorbb}[3]
{\mbox{$#1 \rightarrow #2 \rightarrow #3$}}

\newcommand{\threehorbs}[4]
{\mbox{$#1 \rightarrow #2 \stackrel{#4}{\rightarrow} #3$}}

\newcommand{\threehorsb}[4]
{\mbox{$#1 \stackrel{#4}{\rightarrow} #2 \rightarrow #3$}}

\newcommand{\threehorss}[5]
{\mbox{$#1 \stackrel{#4}{\rightarrow} #2 \stackrel{#5}{\rightarrow} #3$}}

\newcommand{\threehorll}[5]
{\mbox{$#1 \stackrel{#4}{\longrightarrow} #2
\stackrel{#5}{\longrightarrow} #3$}}

\newcommand{\trianglea}[6]
{\begin{center}
 \begin{picture}(50,48)(-25,-24)
\setlength{\unitlength}{1.0mm}
\put(-18,9){\makebox(0,0){$#1$}}
\put(18,9){\makebox(0,0){$#2$}}
\put(0,-9){\makebox(0,0){$#3$}}
\put(-13,4){\vector(1,-1){8}}
\put(13,4){\vector(-1,-1){8}}
\put(-11,9){\vector(1,0){22}}
\put(0,11){\makebox(0,0)[b]{$#4$}}
\put(-10,-1){\makebox(0,0)[tr]{$#5$}}
\put(10,-1){\makebox(0,0)[tl]{$#6$}}
\end{picture}
\end{center}}

\newcommand{\twohorb}[2]
{\mbox{$#1 \rightarrow #2$}}

\newcommand{\twohors}[3]
{\mbox{$#1 \stackrel{#3}{\rightarrow} #2$}}

\newcommand{\twohorl}[3]
{\mbox{$#1 \stackrel{#3}{\longrightarrow} #2$}}
\newcommand{\N}{\mbox{$I\!\! N$}}	
\newcommand{\Z}{\mbox{$Z\!\!\! Z\!$}}	
\newcommand{\Q}{\mbox{$Q\!\!\!\! I$}}	
\newcommand{\R}{\mbox{$I\!\! R$}}	
\newcommand{\mathR}{R\!\!\!\! I \,\,}

\newcommand{\Proj}{\mbox{$I\!\! P$}}	
\newcommand{\mathC}{C\!\!\!\! I \,\,}
\newcommand{\C}{\mbox{$\, I\!\!\!\! C $}}	

\newcommand{\F}{\mbox{$I\!\! F\!$}}	
\newcommand{\putcc}[3]{\put(#1,#2){\makebox(0,0){$#3$}}}
\newcommand{\putcr}[3]{\put(#1,#2){\makebox(0,0)[r]{$#3$}}}
\newcommand{\putcl}[3]{\put(#1,#2){\makebox(0,0)[l]{$#3$}}}
\newcommand{\puttc}[3]{\put(#1,#2){\makebox(0,0)[t]{$#3$}}}
\newcommand{\puttr}[3]{\put(#1,#2){\makebox(0,0)[tr]{$#3$}}}
\newcommand{\puttl}[3]{\put(#1,#2){\makebox(0,0)[tl]{$#3$}}}
\newcommand{\putbc}[3]{\put(#1,#2){\makebox(0,0)[b]{$#3$}}}
\newcommand{\putbr}[3]{\put(#1,#2){\makebox(0,0)[br]{$#3$}}}
\newcommand{\putbl}[3]{\put(#1,#2){\makebox(0,0)[bl]{$#3$}}}
\newcommand{\la}{\mbox{\bf g}}
\newcommand{\Kahler}{K\"{a}hler}
\newcommand{\calAA}{\mbox{$\cal A$}}
\newcommand{\calE}{\mbox{$\cal E$}}
\newcommand{\calS}{\mbox{$\cal S$}}
\newcommand{\calG}{\mbox{$\cal G$}}
\newcommand{\calK}{\mbox{$\cal K$}}
\newcommand{\calM}{\mbox{$\cal M$}}
\newcommand{\calN}{\mbox{$\cal N$}}
\newcommand{\calL}{\mbox{$\cal L$}}
\newcommand{\calT}{\mbox{$\cal T$}}
\newcommand{\calHH}{\mbox{$\cal H$}}

\newcommand{\calGC}{\mbox{$\cal G$$^{\mathC}$}}
\newcommand{\calNs}{\mbox{$\cal N$$_{S}$}}
\newcommand{\calCs}{\mbox{$\cal C$$_{S}$}}

\newcommand{\calAint}{\mbox{$\cal A$$^{1,1}$}}
\newcommand{\End}{\mbox{\/\/End\/\/}_{0}}
\newcommand{\tcE}{\tilde{\calE}}
\newcommand{\tcF}{\tilde{\calF}}

\newcommand{\tum}{\tilde{U}-V_{U}}
\newcommand{\flag}{\, \mbox{flag} \,}
\newcommand{\delbarF}{\delbar_{\cal{F}}}
\newcommand{\delbarE}{\delbar_{\cal{E}}}

\newtheorem{theo}{Theorem}[section]
\newtheorem{prop}[theo]{Proposition}
\newtheorem{rema}[theo]{Remark}
\newtheorem{defi}[theo]{Definition}
\newtheorem{coro}[theo]{Corollary}
\newtheorem{lemm}[theo]{Lemma}


\newcommand{\Aut}{\mbox{Aut}}
\newcommand{\bC}{C\!\!\!\! \underline{I} \,\,}
\newcommand{\blockdiag}{\mbox{block-diag}\,}
\newcommand{\cc}{\mbox{c}}
\newcommand{\calF}{\mbox{$\cal F$}}
\newcommand{\calO}{\mbox{$\cal O$}}
\newcommand{\calW}{\mbox{$\cal W$}}
\newcommand{\diag}{\mbox{diag}\,}
\newcommand{\Endo}{\mbox{End}}
\newcommand{\Gl}{\mbox{Gl}}
\newcommand{\HH}{\mbox{H}}
\newcommand{\Hol}{\mbox{Hol}}
\newcommand{\Image}{\mbox{Im}\,}
\newcommand{\m}[1]{\mbox{#1}}
\newcommand{\mod}{\,\,\, \mbox{mod}\,\,\,}
\newcommand{\normlog}{\mbox{norm log}\, }
\newcommand{\parab}{\mbox{par}\,}
\newcommand{\para}{\mbox{P}\,}
\newcommand{\PD}{\mbox{PD}\,}
\newcommand{\rank}{\mbox{rank}\,}
\newcommand{\Real}{\mbox{Re}\,}
\newcommand{\slope}{\mbox{slope}\,}
\newcommand{\st}{\mbox{ s.t. }\,}
\newcommand{\SU}{\mbox{SU}\,}
\newcommand{\Tr}{\mbox{Tr\,}}
\newcommand{\Tu}{\tilde{U}}
\newcommand{\zP}{z^{\Phi}}
\newcommand{\zmP}{z^{-\Phi}}

\newcommand{\showlabel}[1]{ 
 \label{#1}}
\newcommand{\dimension}{\mbox{dim}}
\newcommand{\stopremark}{ \end{rema} }
\newcommand{\remark}{ \begin{rema} \em  }
\newcommand{\proof}{\noindent{\em Proof:} }
\newcommand{\proofof}[1]{{\noindent \bf Proof} (of #1){\bf :} }
\newcommand{\stopproof}{$\Box$ \hfill\vspace{0.5cm} }
\newcommand{\tu}{\mbox{$\cal{T}(\cal{U})$}}
\newcommand{\ms}{\mbox{$\cal{M}_{S}$}}
\newcommand{\ps}{\mbox{$\cal{P}_{S}$}}
\newcommand{\es}{\mbox{$\cal{E}_{S}$}}
\newcommand{\vv}{\mbox{$\/ \cal V$}}
\newcommand{\uu}{\mbox{$\/ \cal U$}}
\newcommand{\tilu}{\mbox{$\tilde{U}$}}
\newcommand{\M}{\mbox{$\Sigma$}}
\newcommand{\G}{\mbox{$\Gamma$}}
\newcommand{\D}{\mbox{$\Delta$}}
\newcommand{\h}{\mbox{$\lambda$}}
\newcommand{\g}{\mbox{$\gamma$}}
\newcommand{\Gx}{\mbox{$\Gamma_{x}$}}
\newcommand{\GGx}{\mbox{$G_{x}$}}
\newcommand{\calH}{\mbox{$\cal H$}}

\newcommand{\Gy}{\mbox{$\Gamma_{y}$}}
\newcommand{\Gz}{\mbox{$\Gamma_{z}$}}
\newcommand{\sx}{\mbox{$S_{x}$}}
\newcommand{\sy}{\mbox{$S_{y}$}}
\newcommand{\Sx}{\mbox{$\cal S$$_{x}$}}
\newcommand{\Sy}{\mbox{$\cal S$$_{y}$}}
\newcommand{\Se}{\mbox{$\cal S$$_{0}$}}

\newcommand{\al}{\mbox{$\alpha$}}
\newcommand{\shortfunction}[5]{\mbox{$ #1:#2
\rightarrow #3:#4 \mapsto #5 $} }
\newcommand{\function}[5]{\[ \begin{array}{ccccl}
                              #1 & : & #2 & \rightarrow & #3 \\[1.5mm]
                                 &   & #4 & \mapsto     & #5
                             \end{array}
                          \]}
\newcommand{\action}[4]{\[ \begin{array}{ccl}
                               #1 & \rightarrow & #2 \\[1.5mm]
                               #3 & \mapsto     & #4
                             \end{array}
                          \]}

\newcommand{\del}{\mbox{$\partial$}}
\newcommand{\dd}{\mbox{d}}
\newcommand{\delbar}{\mbox{$\bar{\partial}$}}
\newcommand{\mathdelbar}{\bar{\partial}}
\newcommand{\homo}{\mbox{Hom}}
\newcommand{\olog}{\Omega^{1}( \mbox{log} \,\, 0)}
\newcommand{\logd}{\frac{\dd z}{z}}
\newcommand{\diagmatrix}[2]{\left( \begin{array}{ccccc}
                                    #1 & & & & \\
                                       &\cdot & &0 & \\
                                       & & \cdot & & \\
                                        &0 & & \cdot & \\
                                       & & & & #2
                                   \end{array}
                            \right)
                           }
\newcommand{\uppermatrix}[2]{\left( \begin{array}{ccccc}
                                    #1 & & & & \\
                                       &\cdot & & \star  & \\
                                       & & \cdot & & \\
                                        & 0 & & \cdot & \\
                                       & & & & #2
                                   \end{array}
                            \right)
                           }
\newcommand{\twomatrix}[4]{\left( \begin{array}{cc}
                                   #1 & #2 \\
                                    #3 & #4
                                  \end{array}
                            \right)
                            }
\newcommand{\blanctwomatrix}[4]{ \begin{array}{cc}
                                   #1 & #2 \\
                                    #3 & #4
                                  \end{array}
                            }

\newcommand{\bigmatrix}[4]{\left( \begin{array}{cccc}
                                  #1 & \cdots & \cdots & #2 \\
                                  \vdots & \ddots & & \vdots \\
                                  \vdots & & \ddots & \vdots \\
                                   #3 & \cdots & \cdots & #4
                                   \end{array}
                           \right)
                             }
\newcommand{\twovector}[2]{\left( \begin{array}{c}
                                      #1 \\
                                      #2
                                \end{array}
                           \right)
                            }
\newcommand{\threematrix}[9]{\left( \begin{array}{ccc}
                                    #1 & #2 & #3 \\
                                    #4 & #5 & #6 \\
                                    #7 & #8 & #9
                                     \end{array}
                               \right)
                             }

\newcommand{\tz}{\tilde{z}}


\noindent{\em Mathematics Subject Classification
(1991) 14J27 32L07 14H60 14D20 }

\footnotetext[1]{The first author was fully supported by the
{\em Rhodes Trust}, Oxford, and the
         {\em Alfried Krupp von Bohlen und Halbach -- Stiftung}, Essen}

\section{Introduction}

If $q:Y \rightarrow \Sigma$ is an elliptic surface (to be made precise)
then the induced map of fundamental groups is an isomorphism if we consider
$\Sigma$ as an orbifold,  \cite{ue}, \cite{dol}.
Hence, we obtain a correspondence
of flat bundles (by bundles we always mean complex vector bundles).
Donaldson showed that each stable degree zero bundle over
$Y$ or $\Sigma$ admits  a unique
Hermitian-Yang-Mills (or anti-self-dual) connection.
In fact, the obvious conditions on a bundle $E' \rightarrow Y$ to come
from $\Sigma$, namely $\cc_{1}(E') \in q^{*}\HH^{2}(\Sigma)$ and
$\cc_{2}(E')=0$, imply that this H.Y.M. connection
is flat. So, there is a correspondence of stable degree zero
bundles over $Y$ and $\Sigma$.
The generalisation to all degrees
 has been shown by Bauer, \cite{bau}, via a direct proof in
algebraic geometry.

Donaldson's result has been
extended by several authors to parabolic bundles,
theorem \ref{b}. This and the use of
extension results for flat bundles, theorem \ref{y},
 gives an identification of stable parabolic
degree zero bundles over $(\Sigma,P)$ and $(Y,P')$,
where $P$ is a finite collection of generic points
in $\Sigma$ and $P'=q^{-1}(P)$,
with representations of the fundamental groups of
$\Sigma-P$ and $Y-P'$, respectively.
These groups are again isomorphic and so,
we have a correspondence of stable parabolic degree
zero bundles over $(Y,P')$ and over $(\Sigma,P)$.
This extends readily to all degrees by tensoring with
 a parabolic line bundle.

Another proof for genuine bundles of any degree,
not using parabolic bundles at all,
relies on the correspondence between
stable bundles over $\Sigma$ ($Y$) and
representations of the fundamental group of circle
bundles (i.e. Seifert fibred spaces) over $\Sigma$ ($Y$), \cite[p 63]{fas},
\cite{bao}. These groups are also isomorphic,
see the proof of proposition \ref{v}.

Combining Bauer's arguments with Kronheimer \&
Mrowka's description of the moduli
spaces of stable parabolic bundles, we show
that these are complex manifolds if we fix
determinants and if they are pull backs from $\Sigma$.
Finally, we consider smooth parabolic bundles and produce some
details about the correspondence of stable ones.

After recalling elliptic surfaces and parabolic bundles in
 the following two sections
we take one section to state the results.
More details of our work can be found in \cite{gan}.

\hfill

\noindent
{\bf Acknowledgements}

\vspace{0.2cm}

We are grateful to P. Kronheimer,
S. Bauer, M. L\"ubke, T. Peternell,
C. Okonek, S. Agnihotri and R. Plantiko for helpful remarks.

\section{Elliptic surfaces}

Throughout, let $q:Y \rightarrow \Sigma$ be an elliptic surface,
i.e. $Y$ is a compact complex
surface, $\Sigma$ a compact Riemann
surface and $q^{-1}(\sigma)$ an elliptic curve for generic,
i.e. all but finitely many, $\sigma \in \Sigma$, c.f.
\cite{gri}, \cite{bpv}.
We will assume that any non-generic
fibre is either a rational curve of multiplicity one
with one  self-intersection (called singular fibre)
 or  a multiple elliptic curve and furthermore, that
there is at least one singular fibre.
Moishezon shows that all elliptic surfaces are deformation
equivalent to the ones we consider, \cite{moi}.

\begin{theo}[Ue, Dolgachev]
\showlabel{k}
If $q:Y \rightarrow \Sigma$ has singular fibres
and if $U_{0} \subseteq \Sigma$ is an open ball such that
$\pi^{-1}(U_{0})$ contains all singular
fibres but no multiple elliptic curves then
$\pi^{-1}(U_{0})$ is simply connected.
\end{theo}

If $Y_{\sigma}:=q^{-1}(\sigma)$ has
multiplicity $m>1$, let $\Tu$ and $B$ be open discs in $\C$,
$\phi:B \rightarrow U$ a chart with $U \subseteq \Sigma$ and
$\phi(0)=\sigma$ and construct a uniformization of $U$ by
\[ \threehorss{\Tu}{B}{U}{z^{m}}{\phi} \]
where $\langle \eta=e^{2 \pi i/m} \rangle = \Z_{m} \subseteq \C$
acts on $\Tu$ in the standard way.
After this construction is done for
all multiple points $\sigma \in \Sigma$ we
think of $\Sigma$ as an orbifold, cf. \cite{fas}.
Bauer points out that for all orbifold Riemann surfaces there exists
an elliptic surface over it.

The natural gauge-theoretic objects on orbifolds are V-bundles:
a (local, complex) rank $r$ V-bundle $E|_{U}$ is isomorphic to
\[ (\Tu \times \C^{r},\Z_{m})
\rightarrow (\Tu,\Z_{m}) \,\,\,\,\,\,\, \mbox{with}\]
 \[ \eta (\tilde{u},z_{1},...,z_{r})=( \eta \tilde{u},
\eta^{a_{1}} z_{1},...,
       \eta^{ a_{r}} z_{r}) \]
for some isotropies $(a_{1},...,a_{r}) \in \{ 0,...,m-1 \}^{r}$.

\begin{theo}[Furuta \& Steer, Seifert]
\showlabel{i}
Smooth V-bundles over $\Sigma$ are classified by rank,
degree (which is rational) and isotropies.
\end{theo}

For any $y \in Y_{\sigma}$ we can choose coordinates
$(z_{1},z_{2})$ on $U' \ni y$ such that
$\phi^{-1} \circ q (z_{1},z_{2})=z^{m}_{2}$.
Hence we can lift $\phi^{-1} \circ q$ locally to a regular
map \[ \sqrt[m]{\phi^{-1} \circ q}=z_{2}:U' \rightarrow \Tu \]
uniquely up to the action of $\Z_{m}$.
In particular, $q$ is a map of orbifolds.

A divisor on $\Sigma$ can be represented by a finite sum
$D=\sum_{i \in I} \sigma_{i}n_{i}/m_{i}$ where
$n_{i} \in \Z$ and $m_{i}$ the multiplicity
of $\sigma_{i} \in \Sigma$.
The vertical divisors on $Y$ are precisely the
pull backs of divisors on $\Sigma$.
Hence, \cite{nas}, the line bundles $\calO (D')$
over $Y$ with vertical divisor
correspond to the line V-bundles over $\Sigma$.

As $q$ is regular away from finitely many points theorem
\ref{k} implies that
$q_{*}:\pi_{1}(Y) \rightarrow \pi_{1}^{V}(\Sigma)$ is an isomorphism.
Here, the orbifold fundamental  group $\pi_{1}^{V}(\Sigma)$ is
an extension of $\pi_{1}(\Sigma)$ by a
torsion  group with
unipotent elements corresponding to the multiple points.
The isomorphism implies that
 the first betti number of $Y$ is even and therefore
$Y$ is Kahler, by \cite{miy}.
 The correspondence between
representations of the fundamental group and flat bundles on manifolds
extends to orbifolds via the construction
with simply connected coverings.
Hence, the flat bundles over $Y$
correspond to the flat V-bundles over $\Sigma$.

\section{Parabolic bundles}

Let $P:=\{ p_{1},...,p_{n} \}$ be a collection of
generic points in $\Sigma$ and put $P':=
\{ P_{j}':=q^{-1}(p_{j}) \}_{j=1...n} \subseteq Y$.
A parabolic V-bundle $E$ over the
pair $(\Sigma,P)$ is a V-bundle $E \rightarrow \Sigma$ with
proper filtrations
\[ E|_{P_{j}}=E_{j,1} \supset E_{j,2} \supset ...
\supset E_{j,l_{j}} \neq 0 \]
and weights
\[ 0 \leq \alpha_{j,1} < \alpha_{j,2} < ... < \alpha_{j,l_{j}} <1 \]
for all $j=1,...,n$.
We call $\mu_{j,k}:=\rank (E_{j,k}/E_{j,k+1})$ the
multiplicity of $\alpha_{j,k}$. Let $\alpha_{j}$ be
the diagonal matrix of rank equal to $\rank E$
and with entries $\alpha_{j,k}$ ($k=1,...,l_{j}$)
 with multiplicities.
Then put $\alpha(E):=\{ \alpha_{1},...,\alpha_{n} \}$.
We will sometimes write $|E|$ for the underlying genuine V-bundle.

Correspondingly, one defines parabolic
bundles $E'\rightarrow (Y,P')$ where the $\{E'_{j,k} \}$
are bundles over $P_{j}'$ with
weights $\{ \alpha'_{j,k}\}$ of  multiplicities
$\{ \mu'_{j,k} \}$, encoded in $\alpha'(E')$.
If the bundle and the filtrations are holomorphic we
speak of a holomorphic parabolic bundle.
Thinking of vectors in $E|_{P}$ as
having (the obvious) weights ($+\infty$ being the weight of
zero vectors), a morphism of parabolic bundles
is a bundle map not decreasing the weight of any
vector. In a direct sum of parabolic bundles the
weight of a vector is the minimum of the weights of it's
projections.

Let us fix a Kahler metric on $Y$
with (1,1)-form $\omega'$. When we use ordinary bundle invariants
and operators (like $\cc_{i}$, $\det$ or $\deg$) on parabolic bundles
we mean to apply them to the underlying bundles.
Then we have
\[ \parab \cc_{1} (E):= \cc_{1}(E) + \sum_{j=1}^{n} \Tr (\alpha_{j})
   \PD (p_{j} ) \in \HH^{2}(\Sigma,\R) \]
and similarly,  $\parab \cc_{1}(E') \in \HH^{2}(Y,\R)$.
Let $\parab \deg E:= \langle \parab \cc_{1}(E), \Sigma \rangle \in \R$ and
$\parab \deg E':=\langle \parab \cc_{1}(E') \cup \omega',Y \rangle \in \R$.
We also have (since $P \cdot P=0$)
\[ \parab \cc_{2}(E'):=\cc_{2}(E') + 2 \sum_{j=1}^{n}
\sum_{k=1}^{l_{j}} \alpha'_{j,k} \PD(d_{j,k}')
    +\frac{1}{2} \parab \cc_{1}^{2}(E') \in \HH^{4}(Y,\R) \]
where $d'_{j,k}:= \deg (E'_{j,k}/E'_{j,k+1})$.

\begin{defi}
A holomorphic parabolic V-bundle $\calE \rightarrow (\Sigma,P)$
is called {\bf stable}
 if for all non-zero
parabolic maps $\calF \rightarrow \calE$,
injective over some point in $\Sigma$ and with
$\rank \calF < \rank \calE$ we have
\[\mbox{deg}\, \calF / \rank \calF < \mbox{deg}\, \calE / \rank \calE. \]
Similarly for $\calE' \rightarrow (Y,P')$.
\end{defi}

\section{Statement of results}

Being a map of orbifolds,
$q:Y \rightarrow \Sigma$ induces pull backs of
holomorphic parabolic V-bundles over
$\Sigma$ to holomorphic parabolic bundles over $Y$.
Our main result is:

\begin{theo}
\showlabel{z}
Pulling back induces a correspondence between  stable
parabo\-lic V-bundles $\calE \rightarrow (\Sigma,P)$
and those stable parabolic bundles $\calE' \rightarrow (Y,P')$
with $\parab \cc_{2}(\calE')=0$ and
$\parab \cc_{1}(\calE') \in q^{*} \HH^{2}(\Sigma,\R)$.
\end{theo}

In particular:

\begin{theo}
\showlabel{a}
Pulling back induces a bijection between
 stable V-bundles $\calE$ over $\Sigma$ and those stable
bundles $\calE'$ over $Y$ with $\cc_{2}(\calE')=0$
and $\det \calE' = \calO (D')$ for a vertical divisor $D'$.
\end{theo}

Theorem \ref{a} has been shown
(under some assumptions on the Kahler metric of $Y$,
\cite[p 511]{bau}) by Bauer first,
using algebraic geometry. Our goal is a differential-geometric proof of
theorem \ref{z} using the correspondence of stable parabolic
bundles and parabolic H.Y.M. connections.

Let $\calE \rightarrow (\Sigma,P)$ be a
stable parabolic bundle, $\calE':=q^{*} \calE$,
$\calT:=\,\, \mbox{Par End}_{0} \calE$
and $\calT':=q^{*} \calT = \,\, \mbox{Par End}_{0} \calE'$.
The important thing to note is that $\calT$
and $\calT'$ are holomorphic (V-)bundles.
We have the deformation complexes (of smooth sections)

\begin{center}
\setlength{\unitlength}{1.5mm}
\begin{picture}(0,30)(0,-15)
\putcc{-26}{10}{\Omega_{\Sigma}^{0} (\calT)}
\putcc{0}{10}{\Omega_{\Sigma}^{0,1} (\calT)}
\putcc{26}{10}{0}
\putcc{-26}{-10}{\Omega_{Y}^{0} (\calT')}
\putcc{0}{-10}{\Omega_{Y}^{0,1} (\calT')}
\putcc{26}{-10}{\Omega_{Y}^{0,2} (\calT')}
\put(-26,6){\vector(0,-1){12}}
\put(0,6){\vector(0,-1){12}}
\put(26,6){\vector(0,-1){12}}
\put(-18,10){\vector(1,0){9}}
\putbc{-15}{12}{\delbar_{\cal{T}}}
\put(-18,-10){\vector(1,0){8}}
\putbc{-13}{-8}{\delbar_{\cal{T}'}}
\put(9,10){\vector(1,0){10}}
\putbc{13}{12}{\delbar_{\cal{T}}}
\put(8,-10){\vector(1,0){9}}
\putbc{13}{-8}{\delbar_{\cal{T}'}}
\putcl{2}{0}{q^{*}}

\end{picture}
\end{center}

Let $E:=_{C^{\infty}} \calE$, $\calO(D)=\det \calE$
and $\calM(E,D)$ be the space of stable
structures on $E$ with determinant $\calO(D)$.
Similarly, define $\calM(E', D')$ for $E':=q^{*}E$
and $D':=q^{*}D$.
Using the extension of standard deformation theory,
\cite{dak}, to parabolic bundles makes
$\calM(E',D')$ into a Hausdorff complex space,
see \cite[p 83]{km2} and  \cite{mun}, with a description
near $\calE'$ given by the zero set of a holomorphic map
\[ \HH^{0,1}_{\mathdelbar}(Y,\calT')
\rightarrow \HH^{0,2}_{\mathdelbar}(Y,\calT'). \]
By the vanishing of the second cohomology
over $\Sigma$, $\calM(E,D)$ is even a complex manifold
with a chart near $\calE$ given by $\HH^{0,1}_{\mathdelbar}(\Sigma,\calT)$.

\begin{prop}[Bauer]
\showlabel{t}
Stability of $\calE$ implies that
\[ q^{*}:\HH^{0,1}_{\mathdelbar}(\Sigma,\calT)
\rightarrow \HH^{0,1}_{\mathdelbar}(Y,\calT') \]
is an isomorphism.
\end{prop}

This proposition and theorem \ref{z} imply:

\begin{coro}
If $\calM_{\cal{E}}(E,D)$ is the connected component
of $\calE$, correspondingly for $\calE'$,
then \[ q^{*}:\calM_{\cal{E}}(E,D) \rightarrow \calM_{\cal{E}'}(E',D') \]
is an isomorphism of complex manifolds.
\end{coro}

\begin{prop}
\showlabel{u}
If $E' \rightarrow (Y,P')$ is a rank $r$ parabolic bundle satisfying
$\parab \cc_{2} (E')=0$, $\parab \cc_{1} (E')
\in q^{*} \HH^{2}(\Sigma,\R)$ and if the space
$\calM(E')$ of stable structures on $E'$ is non-empty then
\begin{description}
\item[(i)] there exists a unique
line V-bundle $L \rightarrow \Sigma$ with $q^{*} L=\det E'$.
       Also, $E'$ is uniquely determined by $L$ and it's weights $\alpha'$;
\item[(ii)] we have
\[ \calM(E') = \bigsqcup_{a} \calM (L,r,a,\alpha') \]
where the union is over
isotropies $a \in \Z_{m_{1}} \times ... \times \Z_{m_{k}}$
($k=\sharp \{ $ marked points on $\Sigma \}$) compatible with $L$ and
$\calM(L,r,a,\alpha')$ is the space of stable
structures on the unique rank $r$ parabolic V-bundle over
$(\Sigma,P)$ with determinant $L$, isotropies $a$ and weights $\alpha'$.
\end{description}
\end{prop}

\section{Stable parabolic bundles}

\begin{defi}
A {\bf (parabolic) hermitian metric} $h'$ on
a parabolic bundle $E'\rightarrow (Y,P')$
(similarly for $E\rightarrow (\Sigma,P)$) is
a hermitian metric on $E'|_{Y-P'}$ such that
\[ \forall y \in P'_{j} \,\,\,\,\,\,\,\,
   E'_{j,k}(y)=\{ s(y) \,\, | \,\, s
\in \Gamma_{loc}(E') \st h'(s(-))=O(\dd (P'_{j},-)^{\alpha'_{j,k}})
  \}. \]
\end{defi}

For a holomorphic parabolic bundle $\calE'$, $h'$ induces
a Chern connection on $\calE'|_{Y=P'}$.
So, one can talk of  H.Y.M.-connections.
Around $P'_{j}$, a parabolic connection has holonomy
conjugated to $\alpha'_{j}$, \cite{kro1},
and can  therefore not be extended over $P'_{j}$, in general.
A parabolic connection is called reducible if the
parabolic bundle decomposes together with the
connection.
The bridge between algebraic and differential geometry is

\begin{theo}[\cite{mas}, \cite{saw}, \cite{mun}, \cite{b93}, \cite{nas}]
\showlabel{b}
Any degree zero
holomorphic parabolic bundle over a
complex Kahler surface or orbifold of dimension one
is stable if and only if it admits an
irreducible H.Y.M. metric, unique up to isomorphism.
Furthermore, an H.Y.M. connection has finite action.
\end{theo}

The primary result, for genuine bundles, is due to Donaldson;
 Munari's proof relies on Simpson's
work.
A H.Y.M. connection over $\Sigma$ is obviously flat.

\begin{lemm}
\showlabel{c}
If $\calE' \rightarrow (Y,P')$ has
$\parab \cc_{2}(\calE')=0$ and $\parab \cc_{1}(\calE') \in
q^{*} \HH^{2}(\Sigma,\R)$ then a H.Y.M. connection is neccessarily flat.
\end{lemm}

\proof
By \cite[p 100]{mun} or \cite{gan}
 we obtain the second (in particular the existence of) and, from \cite{dak},
the third equality in
\[0=8 \pi^{2} (\parab \cc_{2}(\calE') - \frac{1}{2}
\parab \cc_{1}^{2}(\calE'))=
 \int_{Y-P'} \Tr (F^{2}) = ||F^{-}||^{2} - ||F^{+}||^{2}. \]
\stopproof

\begin{theo}[Munari, Biquard, Simpson]
\showlabel{y}
A holomorphic bundle over
$Y-P'$ ($\Sigma-P$) with a flat hermitian metric  extends uniquely
(up to isomorphism) to a holomorphic
parabolic bundle over $(Y,P)$ ($(\Sigma,P)$) such that
the hermitian metric becomes a parabolic metric.
\end{theo}

Now we prove theorem \ref{z}:

\proof
At first, we treat the special case of degree zero bundles.
If $\calE \rightarrow (\Sigma,P)$ is stable and of
degree zero, let $A$ be the unique irreducible flat parabolic
connection on $\calE$. By the uniqueness
in theorem \ref{y}, $A|_{\Sigma-P}$ is still irreducible.
The regularity of $q|_{Y-P'}$ away from finitely many
points and theorem \ref{k}
imply that $(q|_{Y-P'})_{*}:\pi_{1}(Y-P')
\rightarrow \pi_{1}^{V}(\Sigma-P)$
is an isomorphism.
Hence, if $\calE':=q^{*} \calE$ and $A':=q^{*}A$
then $A'$ is an irreducible flat parabolic connection
and hence $\calE'$ is stable.

Conversely, if $\calE' \rightarrow (Y,P')$
is stable and of degree zero, let $A'$ be the unique
parabolic H.Y.M. connection, which is
flat by lemma \ref{c} and $A'|_{Y-P'}$ is irreducible.
We push forward
$(\calE',A')|_{Y-P'}$ to $\Sigma-P$.
By theorem \ref{y}, this push forward extends uniquely
to a holomorphic parabolic bundle with irreducible parabolic  connection
$(\calE,A) \rightarrow (\Sigma,P)$.
As $q^{*}(\calE,A)$ and $(\calE',A')$
are isomorphic over $Y-P'$ the uniqueness of theorem \ref{y}
implies that they are isomorphic over $(Y,P')$.
We are finished with the degree zero case.

To show the general case, fix a generic point $p \in \Sigma-P$ and let
$\Delta:=\int_{q^{-1}(p)} \omega'$.
For each $d \in \R$ let $[d]$ be it's
integer part, let $[p]$ be the holomorphic line bundle with
divisor $p$ and define the holomorphic parabolic line bundle
$\calL_{d} \rightarrow (\Sigma,p)$
to be $\calL_{d}:=[p]^{[d]}$ with weight $d-[d]$ over $p$.
Let $\calL'_{d}:=q^{*} \calL_{d}$ be the
parabolic pull back. We have $\parab \deg \calL_{d}=d$ and
$\parab \deg \calL'_{d}=\Delta d$.

Tensoring any rank $r$ stable
parabolic bundle $\calE \rightarrow (\Sigma,P)$ ($\calE' \rightarrow
(Y,P')$) of degree
$-dr$ ($-\Delta d r$) with $\calL_{d}$ ($\calL'_{d}$) gives a stable
parabolic bundle of
degree zero over
$(\Sigma,P \cup \{ p\} )$ ($(Y,P' \cup \{ q^{-1}(p) \} )$).
Furthermore, we have $\parab \cc_{2}(\calE' \otimes \calL'_{d})=0$.
So we are reduced to the degree zero case.
\stopproof

We use Bauer's arguments, \cite[p 514]{bau}, to prove proposition \ref{t}:

\proof
The sheaf $\calT$ can be considered as a sheaf
on $\Sigma$ or on $|\Sigma|$, the underlying Riemann
surface of the orbifold $\Sigma$, see the correspondence
between V-bundles and parabolic
bundles in section 5  of \cite{fas}.
As sheaves, $q_{*} \calT' = \calT$.
The Leray spectral sequence induces an exact sequence, \cite[p 10]{bpv},
\[ 0 \rightarrow \threehorsb{\HH^{1}(\Sigma,\calT)}%
{\HH^{1}(Y,\calT')}{\HH^{0}(\Sigma,q_{*1}\calT')}{q^{*}}
     \rightarrow ... \]
where $q_{*1} \calT'$ is the first direct image sheaf of $\calT'$.
It suffices to see that $\HH^{0}(\Sigma,q_{*1}\calT')=0$.

Relative duality, \cite[p 99]{bpv}, gives
\[ q_{*1}\calT' = q_{*}(\calT'^{*}
\otimes \calK_{Y} \otimes q^{*}\calK^{*}_{|\Sigma |})^{*}=
   \calT \otimes q_{*}(\calK_{Y}
\otimes q^{*} \calK^{*}_{|\Sigma |})^{*} \]
where $\calK_{|\Sigma |}$ is the canonical bundle.
By \cite[p 98, 161-162]{bpv},
\[ \calK_{Y} \otimes q^{*} \calK^{*}_{|\Sigma |} =
q^{*} q_{*1} \calO_{Y}^{*} \otimes
    \calO_{Y}(\sum (m_{i}-1)Y_{\sigma_{i}}) \]
where the sum is over the singular points
$\sigma_{i}$ of $\Sigma$, $q_{*1}\calO_{Y}$ is locally free of rank one
since
all the other sheaves in this identity are, and
 \[ \deg (q_{*1} \calO_{Y})=- \chi (\calO_{Y}).  \]

In particular, $\calK_{Y} \cdot \calK_{Y}=0$.
Hence, $\chi (Y)=12 \chi(\calO_{Y})$, \cite[p 472]{gri},
which is equal to the positive number of singular fibres, cf. \cite{ue}.
Now, $\calO_{Y}((m_{i}-1)Y_{\sigma_{i}})=
q^{*}\calO_{\Sigma}(\frac{m_{i}-1}{m_{i}}\sigma_{i})$.
(This  is in fact a trivial sheaf over $\Sigma$.)
We obtain
\[ q_{*1}\calT'=\calT \otimes q_{*1} \calO_{Y}
\otimes  \calO_{\Sigma}(\sum \frac{1-m_{i}}{m_{i}} \sigma_{i})\]
and any non zero section of this induces a non zero map
$\calE \rightarrow \calE \otimes \calL$
for a negative line V-bundle $\calL$.
This is ruled out by stability of $\calE$.
\stopproof

\section{Smooth parabolic bundles}

\begin{prop}
\showlabel{v}
Two smooth line V-bundles
over $\Sigma$ are isomorphic if their pull backs to $Y$ are isomorphic.
\end{prop}

There is an equivalence between smooth line V-bundles
over $\Sigma$ and $\HH^{2}_{V}(\Sigma,\Z)$, \cite{fas}.
As we don't have a sufficient theory of V-cohomology however,
our proof is not by showing
injectivity of $\HH_{V}^{2}(\Sigma,\Z) \rightarrow \HH^{2}(Y,\Z)$.

\proof
It suffices to show that a smooth line V-bundle
$L \rightarrow \Sigma$ is trivial if $L':=q^{*}L \rightarrow Y$ is trivial.
Let us write $q':SL' \rightarrow SL$ for the induced map of circle bundles.
We use the isomorphism
$q_{*}:\pi_{1}(Y) \rightarrow \pi_{1}^{V}(\Sigma)$ and
the commutative diagram

\setlength{\unitlength}{1.0mm}
\begin{center}
\begin{picture}(0,40)(0,-20)
\putcc{-52}{13}{0}
\putcc{-26}{13}{K'}
\putcc{0}{13}{\pi_{1}(SL')}
\putcc{26}{13}{\pi_{1}(Y)}
\putcc{52}{13}{0}

\putcc{-52}{-13}{0}
\putcc{-26}{-13}{K}
\putcc{0}{-13}{\pi_{1}^{V}(SL)}
\putcc{26}{-13}{\pi_{1}^{V}(\Sigma)}
\putcc{52}{-13}{0}

\put(-45,13){\vector(1,0){12}}
\put(-19,13){\vector(1,0){10}}
\put(9,13){\vector(1,0){10}}
\put(33,13){\vector(1,0){12}}

\put(-45,-13){\vector(1,0){12}}
\put(-19,-13){\vector(1,0){10}}
\put(9,-13){\vector(1,0){10}}
\put(33,-13){\vector(1,0){12}}

\put(-52,8){\vector(0,-1){16}}
\put(-26,8){\vector(0,-1){16}}
\put(0,8){\vector(0,-1){16}}
\put(52,8){\vector(0,-1){16}}
\put(26,8){\vector(0,-1){16}}

\putbc{13}{15}{\pi'_{*}}
\putbc{13}{-11}{\pi_{*}}
\putcl{2}{0}{q'_{*}}
\putcl{28}{0}{q_{*}}

\end{picture}
\end{center}

where $K=\langle k \rangle$, $K'=\langle k' \rangle$
for regular fibres $k$ of $SL$ and $k'$ of $SL'$.
The rows are exact because bundles are always regular maps.
By the five-lemma, $q'_{*}$ is an isomorphism
if $(k' \mapsto k):K' \rightarrow K$ is an
isomorphism. Certainly, it is surjective.
Assume there is a V-homotopy $H: [0,1] \times [0,1] \rightarrow SL$
with boundary $k^{n}$.
Since $q$ is regular away from finitely many points, there exists
$H':[0,1] \times [0,1] \rightarrow Y$ lifting $\pi \circ H$.
Hence there exists
$\tilde{H}:[0,1] \times [0,1] \rightarrow SL'$ lifting $H$ and $H'$.
If $k \in S(L_{x})$ then
$\del H' = \pi' \circ \del \tilde{H} \subseteq Y_{x}$
 and  $\del H =q' \circ \del \tilde{H} =k^{n}$.

W.l.o.g. $x \in U_{0}$, where $U_{0}$ is as in theorem \ref{k},
 and we are working on $SL'|_{q^{-1}(U_{0})}=S^{1} \times q^{-1}(U_{0})$.
So we can lift some homotopy (inside $q^{-1}(U_{0})$)
with boundary $\del H'$ to one relating
$\del \tilde{H}$ to $(k')^{n}$.
Hence, $q'_{*}$ is an isomorphism.

Seifert proved that $\pi_{1}^{V}(SL)=$
\[ \langle a_{j}, b_{j},g_{i},k \,\, : \,\, [a_{j},k]=
[b_{j},k]=[g_{i},k]=1=%
g_{i}^{m_{i}}k^{\beta_{i}}=k^{-b}
\prod_{j=1}^{g} [a_{j},b_{j}] \prod_{i=1}^{n} g_{i} \rangle \]
where $g$ is the genus of $\Sigma$ and
$m_{i}$ is the multiplicity of $\sigma_{i}$.
Furthermore, the isotropy
$\beta_{i} \mod m_{i}$  of $L$ at $\sigma_{i}$ and
$\deg L=b+\sum_{1}^{n} \beta_{i}/m_{i}$ are
independent of the choices of lifts $g_{i}, a_{j}$ and $b_{j}$
of the generators of $\pi_{1}^{V} (\Sigma)$.
By the isomorphism of the above extensions
and if $SL'$ is trivial, we can choose lifts  such that
all $\beta_{i}$'s are zero as well as $b$.
By theorem \ref{i}, this implies that $L$ is trivial.
\stopproof

Now we prove proposition \ref{u}:

\proof
Theorem \ref{z} implies the
existence of some $L \rightarrow \Sigma$ with $\det E'=q^{*}L$ and that the
parabolic filtration of $E'$ along $P'$ is
by trivial bundles; in particluar $\cc_{2} (E')=0$.
Proposition \ref{v} gives the uniqueness of $L$.

Now, $|E'|$ is uniquely determined by $\det E'$. Two different
filtrations of $E'|_{P'_{j}}$ by trivial
subbundles are related by a map $P_{j}' \rightarrow \mbox{Sl}(r,\C)$
which can be extended to an isomorphism
of $E'$ being the identity outside a tubular neighbourhood
of $P'_{j}$ since $\mbox{Sl}(r,\C)$
is simply connected. This shows (i).
After the last argument and by theorem \ref{i},
part (ii) follows from proposition \ref{v} and theorem \ref{z}.
\stopproof


\end{document}